\documentclass{amsart}
\usepackage{amsmath}
  \usepackage{paralist}
  \usepackage{graphics} %% add this and next lines if pictures should be in esp format
  \usepackage{epsfig} %For pictures: screened artwork should be set up with an 85 or 100 line screen
\usepackage{caption}
\usepackage{subcaption}
\usepackage{graphicx}  \usepackage{epstopdf}%This is to transfer .eps figure to .pdf figure; please compile your paper using PDFLeTex or PDFTeXify.
 \usepackage[colorlinks=true]{hyperref}
\hypersetup{urlcolor=blue, citecolor=red}

\def\currentvolume{X}
 \def\currentissue{X}
  \def\currentyear{200X}
   \def\currentmonth{XX}
    \def\ppages{X--XX}

\newtheorem{theorem}{Theorem}[section]

\newtheorem{lemma}[theorem]{Lemma}

\theoremstyle{definition}
\newtheorem{definition}[theorem]{Definition}

\DeclareMathOperator*{\argmin}{arg\,min}
\title[Conservative neural network solutions of kinetic equations] %Use the shortened version of the full title
      {Lagrangian dual framework for conservative neural network solutions of kinetic equations}

\author[Hwijae Son and Hyung Ju Hwang]{}


 \email{son9409@kaist.ac.kr}
 \email{hjhwang@postech.ac.kr}

\thanks{$^*$ Corresponding author: Hyung Ju Hwang}

\allowdisplaybreaks

\begin{document}
\maketitle
 
% Enter the first author's name and address:

\centerline{\scshape Hyung Ju Hwang$^*$}
\medskip
{\footnotesize
 % please put the address of the second  and third author
 \centerline{Department of Mathematics}
   \centerline{Pohang University of Science and Technology, Pohang, Republic of Korea}
}
\medskip

\centerline{\scshape Hwijae Son}
\medskip
{\footnotesize
% please put the address of the first author
\centerline{Stochastic Analysis and Application Research Center}
 \centerline{Korea Advanced Institute of Science and Technology, Daejeon, Republic of Korea}
} % Do not forget to end the {\footnotesize by the sign }

\bigskip

% The name of the associate editor will be entered by an editorial staff
% "Communicated by the associate editor name" is not needed for special issue.
% \centerline{(Communicated by the associate editor name)}

%The abstract of your paper
\begin{abstract}
In this paper, we propose a novel conservative formulation for solving kinetic equations via neural networks. More precisely, we formulate the learning problem as a constrained optimization problem with constraints that represent the physical conservation laws. The constraints are relaxed toward the residual loss function by the Lagrangian duality. By imposing physical conservation properties of the solution as constraints of the learning problem, we demonstrate far more accurate approximations of the solutions in terms of errors and the conservation laws, for the kinetic Fokker-Planck equation and the homogeneous Boltzmann equation.
\end{abstract}

%The title of your section 1
\section{Introduction}
Deep neural networks are ubiquitous tools for many fields of science and engineering, such as computer vision, and natural language processing, over the last decades. However, there have been fewer interests in the application of deep neural networks in the field of scientific computing, despite its exceptional success. More recently, approximating solutions of partial differential equations (PDEs) using neural networks is widely studied. In order to obtain a neural network solution of a given PDE, one needs a loss function which guarantees that minimizing the loss function results in an accurate approximation of an analytic solution. One can think of a residual minimization which minimizes a residual of PDEs \cite{lagaris1998artificial, lagaris2000neural, sirignano2018dgm, raissi2019physics, lu2021deepxde, NHM, hwang2020trend, son2021sobolev}, or a variational form \cite{weinan2018deep, muller2019deep, liao2019deep, kharazmi2021hp}.

Compared to the traditional mesh-based numerical schemes, there are several pros and cons of neural network methods. To state several advantages, neural networks are free of mesh generation by its nature, and hence free of so called curse of dimensionality. We refer to the readers \cite{weinan2017deep, sirignano2018dgm, han2018solving, weinan2018deep} for more information about high-dimensional applications. Another advantage pointed out by \cite{berg2018unified} is that the neural network method is less affected by a domain complexity of the problem. They showed an example of 2-dimensional domain with a highly complicated boundary so that the traditional finite element method fails to apply. However, the authors demonstrated that a neural network can learn a solution of a PDE in that complex domain successfully.

A major disadvantage is that the training of neural network solutions requires a more expensive computational cost compared to the mesh-based schemes. Recently there have been numerous works to improve the convergence of a neural network solution. For instance, \cite{lyu2020enforcing} proposed a way to enforce exact boundary conditions for various kinds of boundary conditions, and claim that the training gets much faster. Another work by \cite{muller2021notes} showed that the convergence mode can be improved from $H^{1/2}$ to $H^2$ when the exact boundary condition is imposed on the Poisson problem. Regarding the convergence of neural network solutions, \cite{wang2020and} theoretically showed that there exists a bias between the residual and the collocation losses, and proposed a method, based on the eigenvalues of its neural tangent kernel, to overcome the training bias. \cite{mcclenny2020self, van2020optimally} also made similar contribution on balancing each of loss components to improve the convergence speed of the neural network solutions.

The aforementioned works to overcome the disadvantages of the neural network methods are more or less generic ones, not equation specific. In this paper, we propose a new class of loss functions which enforce the approximated solution to conserve relevant physical quantities which is an important concept in the kinetic equations. Physical conservation law is one of the fundamental properties of the Boltzmann equation and the other kinetic equations. Therefore, a numerical scheme that satisfies the conservation laws is always of interest. For example, splitting schemes are very popular since they are simple to design to preserve physical quantities (See, \cite{dimarco2014numerical, aristov2012direct} for more information). However, despite its importance, the conservation laws are not yet widely studied in the neural network community. We believe that this is the first attempt to design new loss functions to enforce various kinds of physical conservation laws of kinetic equations.

\subsection{Related works}
Neural network methods for the kinetic equations are only recently developed research areas. The kinetic Fokker--Planck equation was studied in \cite{hwang2020trend} via a neural network approach. The authors showed the existence of a neural network that can arbitrarily minimize the residual loss function and that minimizing the loss function is sufficient to guarantee the convergence of the neural network solution to a priori analytic solution under the specular reflection boundary condition. Moreover, they showed in numerical examples that the neural networks successfully approximate the solutions for various kinds of boundary conditions. There is another work which mainly deals with the kinetic Fokker--Planck equation \cite{son2021sobolev}. The authors proved that the convergence of neural network solutions to a priori analytic solution under the periodic boundary condition and proposed a new loss function that can improve the convergence mode from $L^2$ to $H^1$. Their numerical results show much faster convergence compared to the traditional $L^2$ loss functions. For different type of equations, \cite{lou2020physics} employed physics informed neural network for solving forward and inverse problems related to multi-scale flows with the Bhatnagar-Gross-Krook (BGK) collision model. A recent paper \cite{lee2020model} studies the diffusion limit of the Vlasov--Poisson--Fokker--Planck equation to the Poisson--Nernst--Planck system with simulations via neural networks. 

Regarding the constrained optimization for the neural networks, we want to mention the following studies. There are several works that made their efforts to impose a hard constraint on a neural network, for example, \cite{marquez2017imposing} studied a scalable method for imposing a hard constraint on the output of a neural network. Another work by \cite{ravi2019explicitly} argued that the stochastic gradient descent does not deal with the constraints in a natural way. By incorporating the conditional gradient method, they showed a faster convergence for some specific type of constraints. Another line of researches mainly focus on the Lagrangian duality for constrained optimizations. For example, \cite{fioretto2020lagrangian, nandwani2019primal} reformulated the constrained optimization problem involving neural networks to an unconstrained one by using the Lagrangian relaxation. The changed problem is a minimax problem and thus can be solved by a gradient descent ascent method. Our methodology highly relies on these works. 

\subsection{Outline of the paper}
This paper is organized as follows. In Section \ref{preliminary}, we briefly review the fundamental concept of neural network solutions of PDEs and a general methodology of it. After that, we introduce a well known notion of Lagrangian duality for constrained optimization problems. In Section \ref{methodology} we propose a new class of loss functions and learning framework for training conservative neural network solutions of the kinetic Fokker--Planck equation and the homogeneous Boltzmann equation, which is a main contribution of the paper. In Section \ref{numerical}, we demonstrate superior performances of the proposed methodology in terms of accuracy, and the conservation, through three numerical examples. In Section \ref{discussion} we summarize the results and conclude the paper.

\section{Preliminaries}\label{preliminary}
\subsection{Neural network method for solving PDEs}
We first give a rigorous definition of the fully connected neural network.

\begin{definition}\label{nn_def}
    A fully connected neural network $u_{\theta}(x) : \mathbb{R}^{n_0} \rightarrow \mathbb{R}^{n_L}$ is a function defined by an alternating composition of affine transformations $A_{l}(x) = W_{l}^Tx + b_l$, $W_{l} \in \mathbb{R}^{n_{l-1} \times n_{l}}, b_l \in \mathbb{R}^{n_l}$, and a non-linear activation function $\phi : \mathbb{R} \rightarrow \mathbb{R}$: 
    \begin{equation}
        u_{\theta} (x) = A_L \circ \phi \circ \cdots \circ \phi \circ A_1(x), \nonumber
    \end{equation}
    where $\theta = (vec(W_0), b_0, \cdots ,vec(W_L), b_L)$, and $\phi$ is composed component-wisely.
\end{definition}
Neural networks are known to be universal approximators in corresponding function spaces, and the following theorems summarize the statement.

\begin{theorem}[Theorem 1, in \cite{leshno1993multilayer}]
    Let M denote the set of functions which are in $L_{loc}^{\infty}(\mathbb{R})$ and such that the closure of the discontinuity is of zero Lebesgue measure. Let $\sigma \in M$. Set 
    \begin{equation}
        \Sigma_n = span\{\sigma(w\cdot x + b) : w \in \mathbb{R}^n, b \in \mathbb{R}\}. \nonumber
    \end{equation}
    Then, $\Sigma_n$ is dense in $C(\mathbb{R}^n)$ if and only if $\sigma$ is not an algebraic polynomial.
\end{theorem}
We refer to the readers \cite{cybenko1989approximation, hornik1989multilayer}, for more details about the previous results. Not only for the continuous functions, but neural networks can also approximate arbitrary differentiable function and its derivative simultaneously. The following theorem states the results.

\begin{theorem}[Theorem 2.1, in \cite{li1996simultaneous}]
    Let $K$ be a compact subset of $\mathbb{R}^n$, $n\geq1$, and $f\in \hat{C}^{m_1}(K)\cap \cdots \cap \hat{C}^{m_q}(K)$, where $m_i \in \mathbb{Z}_+^n$ for $1\leq i \leq q$. Also, let $\sigma$ be any non-polynomial function in $C^m(\mathbb{R}^n)$, where $m=\max\{|m_i|, 1\leq i \leq q\}$. Then for any $\epsilon>0$, there is a neural network 
    \begin{equation}
        N(x) = \Sigma_{j=0}^N c_j \sigma(\langle w_j, x\rangle+b_j), \text{ } x\in \mathbb{R}^n, \nonumber
    \end{equation}
    where $c_j \in \mathbb{R}, w_j \in \mathbb{R}^n, b_j \in \mathbb{R}$, such that
    \begin{equation}
        \| D^kf - D^kN \|_{L^{\infty}(K)} \leq \epsilon, \text{ } k \in \mathbb{Z}_+^n, k\leq m_i,\text{ for some }1\leq i\leq q. \nonumber
    \end{equation}
\end{theorem}
Although the above theorems state that neural networks are powerful approximators, how to find such an approximator is not obvious. From now, we introduce a method for solving generic PDEs via neural networks which is an active research area \cite{raissi2019physics, sirignano2018dgm, hwang2020trend, son2021sobolev}.

Consider an initial boundary value problem of a generic PDE :
\begin{equation}
    \begin{split} \label{pde}
        Nu(t,x) &= 0 \text{, for } (t,x,) \in [0,T]\times \Omega, \\
        u(0,x) &= g(x) \text{, for } x \in \Omega \, \\
        Bu(t,x) &= h(t,x) \text{, for } (t,x) \in [0,T]\times\partial\Omega, 
    \end{split}
\end{equation}
where  $N$ is a differential operator, $B$ represents the boundary operator, $g$, and $h$ are initial and boundary conditions, respectively. Starting from the work \cite{lagaris1998artificial}, neural networks are widely considered as trial functions for approximating the solution of \eqref{pde}. 

Let $u_{\theta}(t,x)$ be a neural network with inputs $t$, and $x$ and the parameter vector $\theta = (vec(W_0), b_0, \cdots vec(W_L), b_L)$ as in definition \ref{nn_def}. Define the loss function by penalizing the neural network to satisfy the initial boundary value problem:
\begin{equation}
    Loss(u_\theta) = Loss_{GE}(u_{\theta}) + Loss_{IC}(u_{\theta}) + Loss_{BC}(u_{\theta}), \nonumber
\end{equation}
where 
\begin{equation}
    \begin{split} \label{loss_generic}
        Loss_{GE}(u_{\theta}) &= ||Nu_{\theta}(t,x)||_{L^2([0,T]\times\Omega)}, \\
        Loss_{IC}(u_{\theta}) &= ||u_{\theta}(0,x) - g(x)||_{L^2(\Omega)}, \\
        Loss_{BC}(u_{\theta}) &= ||Bu_{\theta}(t,x) - h(x)||_{L^2([0,T]\times\partial\Omega)}.
    \end{split}
\end{equation}
In practice, we compute \eqref{loss_generic} by its Monte-Carlo approximation. Suppose that we are given uniformly sampled collocation points 
\begin{equation}
    \begin{split}
        \{(t_i, x_j)\}_{i,j=1}^N &\sim UNIF([0,T]\times\Omega), \\
        \{x_j\}_{j=1}^{N_I} &\sim UNIF(\Omega), \nonumber \\
        \{(t_i, x_j)\}_{i,j=1}^{N_B} &\sim UNIF([0,T]\times \partial\Omega),
    \end{split}
\end{equation}
then, we approximate \eqref{loss_generic} by:
\begin{equation}
    \begin{split}\nonumber
        \widehat{Loss}_{GE}(u_\theta) &= \frac{T|\Omega|}{N^2} \sum_{i,j=1}^N |Nu_\theta (t_i,x_j)|^2, \\
        \widehat{Loss}_{IC}(u_\theta) &= \frac{|\Omega|}{N_I} \sum_{j=1}^{N_I} |u_\theta(0,x_j) - g(x_j)|^2, \\
        \widehat{Loss}_{BC}(u_\theta) &= \frac{T|\partial\Omega|}{N_B^2} \sum_{i,j=1}^{N_B} |Bu_\theta(t_i,x_j) - h(t_i, x_j)|^2,
    \end{split}
\end{equation}
and the total loss becomes
\begin{equation}
    \widehat{Loss}(u_\theta) = \widehat{Loss}_{GE}(u_\theta) + \widehat{Loss}_{IC}(u_\theta) + \widehat{Loss}_{BC}(u_\theta).
\end{equation}

After defining the loss function, one need to solve a minimization problem which reads as: 
\begin{equation}
    \argmin_{\theta} \widehat{Loss}(u_\theta). \nonumber
\end{equation}
One can use a gradient-based optimization technique: 
\begin{equation}
    \theta^{new} = \theta^{old} - \eta \nabla_{\theta} \widehat{Loss}(u_\theta) \big\rvert_{\theta=\theta^{old}}, \nonumber
\end{equation}
to minimize the loss function with respect to the neural network parameter $\theta$. In this paper, we fix a first-order gradient-based optimization technique, called ADAM (see, \cite{kingma2014adam} for detailed information). 
\subsection{Constrained Optimization}
Consider a generic constrained optimization problem for an objective function $J(\theta)$, subject to a set of equality constraints $C(\theta) = \{c_1(\theta), \dots, c_n(\theta)\}$ which reads as:
\begin{equation}
    \begin{split} \label{constrained}
        \argmin_\theta \text{ } &J(\theta),  \\
        \text{subject to } &c_1(\theta)= \dots = c_n(\theta)=0. 
    \end{split}
\end{equation}

The simplest and trivial method for solving \eqref{constrained} is to consider a quadratic penalty function. One can reformulate the constrained optimization problem to a unconstrained one by solving:
\begin{equation}\label{penalty}
    \argmin_\theta J(\theta) + \beta^{(k)} \sum_{i=1}^N c_i(\theta)^2. 
\end{equation}
for an increasing scalar sequence $\{\beta^{(k)}\}_{k=1}^{\infty}$ such that $\beta^{(k)} \rightarrow \infty$, as $k \rightarrow \infty$. Although the penalty methods are easy to implement and has advantages of unconstrained nature, it has a severe stability issue due to large values of $\beta^{(k)}$ (See, \cite{bertsekas1976multiplier}).

The method of Lagrange multipliers converts the constrained optimization problem \eqref{constrained} to a unconstrained one by using the Lagrange multiplier. The converted problem can be written as:
\begin{equation}\label{lagrange}
    \argmin_\theta J_\lambda (\theta) = \argmin_\theta J(\theta) + \sum_{i=1}^N\lambda_i c_i(\theta), 
\end{equation}
where $\lambda=(\lambda_1,\dots,\lambda_N) \in \mathbb{R}^N$ is a multiplier vector. One minimizes $J_{\lambda^{(k)}}(\theta)$ for a sequence of multiplier vectors $\{\lambda^{(k)}\}$ generated by: 
\begin{equation}\label{lambda_LM}
    \lambda^{(k+1)}_i = \lambda^{(k)}_i + \eta c_i(\theta^{(k)}), 
\end{equation}
where $\eta$ is a learning rate. Above iteration is a gradient ascent step for the multiplier variable for finding an optimal solution of corresponding dual problem. This method overcomes the stability issue, however, it requires the original problem \eqref{constrained} to have a locally convex structure (See, \cite{luenberger1973introduction, bertsekas1976multiplier}).

Augmented Lagrangian methods combine the ideas of above mentioned methods. In these methods, the penalty term is added to the Lagrangian function $J_\lambda(\theta)$, and forming the augmented Lagrangian function by
\begin{equation}
    J_{\lambda,\beta}(\theta) = J(\theta) + \sum_{i=1}^N \lambda_ic_i(\theta) + \beta\sum_{i=1}^N \|c_i(\theta)\|^2. \nonumber
\end{equation}
An iterative minimization of the form
\begin{equation}
    \argmin_\theta J_{\lambda^{(k)}, \beta^{(k)}}(\theta) = \argmin_\theta J(\theta) + \sum_{i=1}^N \lambda_i^{(k)} c_i(\theta) + \beta^{(k)}\sum_{i=1}^N \|c_i(\theta)\|^2, \nonumber
\end{equation}
is performed and the multiplier sequence $\lambda^{(k)}$ is generated from 
\begin{equation}\label{lambda_AL}
    \lambda^{(k+1)}_i = \lambda^{(k)}_i + \beta^{(k)}c_i(\theta^{(k)}). 
\end{equation}
These type of methods may converge without the need of increasing sequence $\beta^{(k)}$, thus the stability issue can be avoided. Moreover, there is no need for the original problem \eqref{constrained} to have a locally convex structure \cite{bertsekas1976multiplier, sangalli2021constrained}.

\section{Methodology}\label{methodology}
In this section, we first introduce our two model problems and briefly summarize related previous results. Then we introduce our main contribution which makes the training process far more efficient and accurate when solving kinetic equations via neural networks. Our target equations are the kinetic Fokker--Planck equation and the homogeneous Boltzmann equation.

\subsection{Kinetic Fokker--Planck equation} The $d$-dimensional kinetic Fokker--Planck equation reads as:
\begin{equation}\label{fokkerplanck}
    \partial_t f + v\cdot\nabla_x f = \nabla_v \cdot (q \nabla f + p vf), \text{ for } (t,x,v) \in [0,T]\times \Omega \times \mathbb{R}^d, 
\end{equation}
where $\Omega \subset \mathbb{R}^d,$ $q\geq0$ is a diffusion coefficient, $p \geq 0$ is a friction coefficient, and $f=f(t,x,v)$ is a probability density function of particles. In this paper, we consider the $1$-dimensional case in a bounded interval $\Omega=[0,1]$ with the periodic boundary condition. The initial-boundary value problem for the kinetic Fokker--Planck equation reads as:
\begin{equation}
    \begin{split}\label{IBVP_FP}
        \partial_t f + v\partial_xf &= \partial_v(q\partial_v f + p v f), \text{ for } (t,x,v) \in [0,T]\times [0,1] \times \mathbb{R}, \\
        f(0,x,v) &= f_0(x,v) \geq 0, \text{ for } (x,v) \in [0,1]\times\mathbb{R},  \\
        f(t,x,v) &= f(t,1-x,v), \text{ for } (t,x,v) \in \Gamma^-, 
    \end{split}
\end{equation}
where $\Gamma^- = \{(t,0,v) \rvert v>0\} \cup \{(t,1,v) \rvert v<0\}$.
It is well known that the above equation \eqref{fokkerplanck} has a global equilibrium solution called global Maxwellian. Moreover, the conservation of mass and the balance identities for energy and entropy are also given in \cite{soler1997asymptotic}:
\begin{lemma}[Lemma 3.1 in \cite{soler1997asymptotic}]\label{mass_FP}
    Let f be a solution of \eqref{IBVP_FP}, then 
    \begin{equation}
        \frac{d}{dt}\|f(t,\cdot, \cdot)\|_{L^1(\Omega \times \mathbb{R})}=0. \nonumber
    \end{equation}
\end{lemma}

When solving the above initial-boundary value problem \eqref{IBVP_FP} for the kinetic Fokker--Planck equation via neural network, one can consider the loss function as below:
\begin{equation}\label{loss_fp}
    \begin{split}
        Loss_{GE}(f_\theta) &= \|\partial_t f_\theta + v\partial_x f_\theta - \partial_v (qf_\theta+p v f_\theta)\|_{L^2_{t,x,v}}^2, \nonumber\\
        Loss_{IC}(f_\theta) &= \|f_\theta(0,\cdot, \cdot) - f_0(\cdot, \cdot)\|_{L^2_{x,v}}^2, \nonumber\\
        Loss_{BC}(f_\theta) &= \|f_\theta(t,x,v)-f_\theta(t,1-x,v)\|_{L^2(\Gamma^-)},\nonumber
    \end{split}
\end{equation}
where $f_\theta$ denotes a neural network with the parameter vector $\theta$, and the total loss function for the kinetic Fokker--Planck equation:
\begin{equation}\label{total_loss_fp}
    Loss_{FP}(f_\theta) = Loss_{GE}(f_\theta) + Loss_{IC}(f_\theta) + Loss_{BC}(f_\theta).
\end{equation}
This equation is widely studied as a model problem for studying the neural network method for solving kinetic equations. For example, \cite{hwang2020trend} showed the convergence of neural network solutions to the analytic solution when the total loss \eqref{total_loss_fp} vanishes under the specular reflection and the inflow boundary condition, (see, Theorem 3.6 in \cite{hwang2020trend}). Moreover, the authors of \cite{hwang2020trend} demonstrated successful numerical results for different kinds of initial and boundary conditions. More recently, \cite{son2021sobolev} showed the same kind of convergence results for the periodic boundary condition. Furthermore, they showed that the convergence mode can be improved from $L^{\infty}(0,T; L^2_{x,v}))$ to $L^{\infty}(0,T;H^1_{x,v})$ with a slight modification in the loss function (See, Theorem 4.2 in \cite{son2021sobolev}).

The main idea of this paper is to restrict the optimization problem to a smaller function space by using the constrained optimization technique. Thus, we reformulate the problem of finding a neural network solution of \eqref{IBVP_FP} as a constrained optimization problem by:
\begin{equation}
    \begin{split}
        \argmin_\theta\text{ } &Loss_{FP}(f_\theta)\\ \nonumber
        \text{subject to } &\frac{d}{dt} \| f_\theta(t,\cdot, \cdot) \|_{L^1(\Omega \times \mathbb{R})} = 0.
    \end{split}
\end{equation}
We define three loss functions for penalty (P), Lagrange multiplier (L), and augmented Lagrangian method (A) for solving above constrained optimization problem as follows: 
\begin{equation}
    \begin{split}\label{proposed_FP}
        Loss_{FP}^{(P)}(f_\theta) = &Loss_{FP}(f_\theta) + \beta \big\| \frac{d}{dt} \| f_\theta(t,\cdot, \cdot) \|_{L^1(\Omega \times \mathbb{R})} \big\|_{L^2(0,T)}^2,\\
        Loss_{FP}^{(L)}(f_\theta) = &Loss_{FP}(f_\theta) + \lambda(t) \frac{d}{dt} \| f_\theta(t,\cdot, \cdot) \|_{L^1(\Omega \times \mathbb{R})},\\
        Loss_{FP}^{(A)}(f_\theta) = &Loss_{FP}(f_\theta) + \mu \big\| \frac{d}{dt} \| f_\theta(t,\cdot, \cdot) \|_{L^1(\Omega \times \mathbb{R})} \big\|_{L^2(0,T)}^2 \\
        &+ \lambda(t) \frac{d}{dt} \| f_\theta(t,\cdot, \cdot) \|_{L^1(\Omega \times \mathbb{R})},
    \end{split}
\end{equation}
where, $\beta, \mu$ are fixed constants, and $\lambda(t)$ is a multiplier variable which will be updated by gradient ascent as introduced in \eqref{lambda_LM}, and \eqref{lambda_AL}. In practice, we discretize $\lambda(t)$ for fixed grid points $\{t_1, t_2, \dots t_M\}$ and treat it as a vector $(\lambda(t_1), \lambda(t_2), \dots, \lambda(t_M))$. We provide numerical results that show superior performances of the proposed loss function \eqref{proposed_FP} compared to the unconstrained one \eqref{total_loss_fp}. 
\subsection{Homogeneous Boltzmann equation} We also consider the Boltzmann equation in the space homogeneous case
\begin{equation}\label{boltzmann}
    \partial_t f = \frac{1}{\epsilon} Q(f,f), 
\end{equation}
with the initial condition 
\begin{equation}
    f(0, v) = f_0(v), \nonumber
\end{equation}
where $f$ is a non-negative distribution function of particles which move with velocity $v \in \mathbb{R}^3$. The positive constant $\epsilon$ is the Knudsen number and the binary collision operator is given by 
\begin{equation}
    Q(f,f)(v) = \int_{\mathbb{R}^3}\int_{S^2} \sigma(|v-v_{*}|, w)(f(v')f(v'_*) - f(v)f(v_*)) dw dv_*, \nonumber
\end{equation}
where $(v', v'_*)$ is a pair of post-collision velocities which can be written as 
\begin{equation}
    v' = \frac{1}{2} (v+v_* + |v-v_*|w), \quad v'_* = \frac{1}{2} (v+v_* - |v-v_*|w). \nonumber
\end{equation}
The kernel $\sigma$ characterizes the binary interaction of particles. 

During the evoltion process, the collision operator $Q(f,f)$ preserves mass, momentum, and energy %\begin{equation}
%    \begin{split}
%        &\int_{\mathbb{R}^3} Q(f,f) dv = 0, \\
%        &\int_{\mathbb{R}^3} Q(f,f)v dv = 0, \\ \nonumber
%        &\int_{\mathbb{R}^3} Q(f,f)|v|^2 dv = 0,
%    \end{split}
%\end{equation}
and therefore, in the homogeneous case, the following quantities are conserved in time 
\begin{equation}\label{conserve_boltzmann}
    \begin{split}
        &\frac{d}{dt}\int_{\mathbb{R}^3} f dv = 0, \\
        &\frac{d}{dt}\int_{\mathbb{R}^3} fv dv = 0, \\ 
        &\frac{d}{dt}\int_{\mathbb{R}^3} f|v|^2 dv = 0.
    \end{split}
\end{equation}

We define the loss function in the same manner as in \eqref{loss_fp}
\begin{equation}\label{loss_boltzmann}
    \begin{split}
        Loss_{GE}(f_\theta) &=  \|\partial_t f_\theta  - \frac{1}{\epsilon}Q(f_\theta, f_\theta)\|_{L^2_{t,v}}^2 , \\
        Loss_{IC}(f_\theta) &= \|f_\theta(0,\cdot) - f_0(\cdot)\|_{L^2_{v}}^2,
    \end{split}
\end{equation}
and the total loss function for the Boltzmann equation by 
\begin{equation}\label{loss_total_boltzmann}
    Loss_{B}(f_\theta) = Loss_{GE}(f_\theta) + Loss_{IC}(f_\theta).
\end{equation}

We also reformulate the problem as a constrained optimization by employing the identities in \eqref{conserve_boltzmann} as constraints of the optimization problem.
\begin{equation}
    \begin{split}
        \argmin_\theta\text{ } &Loss_{B}(f_\theta)\\ 
        \text{subject to } &\frac{d}{dt} \int_{\mathbb{R}^3} f_\theta(t,v)dv = 0, \\
        &\frac{d}{dt} \int_{\mathbb{R}^3} f_\theta(t,v)vdv = 0, \\
        &\frac{d}{dt} \int_{\mathbb{R}^3} f_\theta(t,v)|v|^2dv = 0, \\
    \end{split}
\end{equation}

We define three loss functions for penalty (P), Lagrange multiplier (L), and augmented Lagrangian method (A) for solving above constrained optimization problem with five constraints as follows: 
\begin{equation}
    \begin{split}\label{proposed_boltzmann}
        Loss_{B}^{(P)}(f_\theta) = &Loss_{B}(f_\theta) + \sum_{i=1}^5 \beta_i \| c_i(t; \theta) \|_{L^2(0,T)}^2, \\
        Loss_{B}^{(L)}(f_\theta) = &Loss_{B}(f_\theta) + \sum_{i=1}^5 \lambda_i(t) c_i(t;\theta),\\
        Loss_{B}^{(A)}(f_\theta) = &Loss_{B}(f_\theta) + \sum_{i=1}^5 \mu \| c_i(t; \theta) \|_{L^2(0,T)}^2 + \sum_{i=1}^5 \lambda_i(t) c_i(t;\theta),
    \end{split}
\end{equation}
where $v=(v_1,v_2,v_3)\in \mathbb{R}^3$, $c_1(t;\theta) = \frac{d}{dt} \int_{\mathbb{R}^3} f_\theta(t,v)dv$, $c_2(t;\theta) = \frac{d}{dt} \int_{\mathbb{R}^3} f_\theta(t,v)v_1dv$, $c_3(t;\theta) = \frac{d}{dt} \int_{\mathbb{R}^3} f_\theta(t,v)v_2dv$, $c_4(t;\theta) = \frac{d}{dt} \int_{\mathbb{R}^3} f_\theta(t,v)v_3dv$, \newline
and $c_5(t;\theta) = \frac{d}{dt} \int_{\mathbb{R}^3} f_\theta(t,v)|v|^2dv$. The parameters $\beta_i, \lambda_i(t)$, and $\mu_i$ are the same as in \eqref{proposed_FP}.

\section{Numerical Results}\label{numerical}
In this section, we demonstrate the superior performance of the proposed constrained optimization method for learning deep neural network solutions of the kinetic PDEs. 
In all three of our numerical results, we use a neural network with 4 hidden layers and 256 nodes for each hidden layer. For the activation function, although we can consider ReLU activation for the Boltzmann equation, but due to the presence of second derivative in \eqref{fokkerplanck} we need at least $C^2$ activation function for the kinetic Fokker--Planck equation. Therefore, we decide to use the hyperbolic tangent function as an activation function. The weights are initialized uniformly, and we use ADAM \cite{kingma2014adam} as an optimizer. For the deep learning framework, we use PyTorch \cite{paszke2019pytorch}. 

\subsection{Kinetic Fokker--Planck equation}
We consider the 1-dimensional kinetic Fokker--Planck equation with the periodic boundary condition and positive initial condition
\begin{equation}
    \begin{split}
        \partial_t f + v\partial_x f &= \partial_v(q\partial_v f + p v f), \text{ for } (t,x,v) \in [0,T]\times [0,1]\times [-V,V], \\\nonumber
        f_(0,x,v) &= f_0(x,v) \geq 0, \text{ for } (x,v)\in [0,1]\times[-V,V],\\
        f(t,x,v) &= f(t,1-x,v), \text{ for } (t,x,v) \in \Gamma^-_V,
    \end{split}
\end{equation}
where we truncate the domain of $v$ variable to $[-V,V]$, and the space $\Gamma^-_V = \{(t,0,v) | 0<v<V\} \cup \{(t,1,v)|-V<v<0\}$ is also truncated accordingly. This truncation of the domain is widely considered in the context of numerical analysis \cite{dimarco2014numerical}, as well as neural network methods \cite{hwang2020trend, son2021sobolev}. In the rest of the paper, we set $V=5$ unless otherwies specified. 

In order to discretize the loss function, we first sample the collocation points from each domain
\begin{equation}
    \begin{split}
        \{(t_i,x_j,v_k)\}_{i,j,k=1}^{N_C} &\sim UNIF([0,T]\times[0,1]\times[-5,5]),\\
        \{(x_j,v_k)\}_{j,k=1}^{N_I} &\sim UNIF([0,1]\times[-5,5]),\\ \nonumber
        \{(t_i,x_j,v_k)\}_{i,j,k=1}^{N_B} &\sim UNIF(\Gamma^-_V),
    \end{split}
\end{equation}
where $N_C, N_I,$ and $N_B$ denote the number of sample points in the whole, initial, and boundary domain, respectively. 

Then we discretize each loss function introduced in section \eqref{total_loss_fp}, and \eqref{proposed_FP} as follows:

\begin{align}
        \widehat{Loss}_{FP} &= \frac{2TV}{N_C^3}\sum_{i,j,k=1}^{N_C}(\partial_t f_\theta + v\partial_x f_\theta - \partial_v(qf_\theta + p v f_\theta))^2\bigg|_{(t_i,x_j,v_k)} \nonumber\\
        &+ \frac{2V}{N_I^2} \sum_{j,k=1}^{N_I}(f_\theta(0, x_j, v_k) - f_0(0,x_j,v_k))^2 \nonumber \\
        &+ \frac{|\Gamma^-_V|}{2N_B^2} \sum_{i,j,k=1}^{N_B} (f_\theta(t_i,x_j,v_k) - f_\theta(t_i,1-x_j,v_k))^2, \nonumber\\
        \widehat{Loss}_{FP}^{(P)} &= \widehat{Loss}_{FP} + \beta \frac{T}{N_C}\sum_{i=1}^{N_C} (\frac{2V}{N_C^2} \sum_{j,k=1}^{N_C}{\frac{d}{dt}f_\theta(t_i, x_j, v_k)})^2, \label{FP_loss}\\
        \widehat{Loss}_{FP}^{(L)} &= \widehat{Loss}_{FP} + \sum_{i=1}^{N_C}\lambda(t_i)  \frac{2V}{N_C^2} \sum_{j,k=1}^{N_C}{\frac{d}{dt}f_\theta(t_i, x_j, v_k)},\nonumber\\
        \widehat{Loss}_{FP}^{(A)} &= \widehat{Loss}_{FP} + \mu \frac{T}{N_C}\sum_{i=1}^{N_C} (\frac{2V}{N_C^2} \sum_{j,k=1}^{N_C}{\frac{d}{dt}f_\theta(t_i, x_j, v_k)})^2 \nonumber\\
        &+ \sum_{i=1}^{N_C}\lambda(t_i)  \frac{2V}{N_C^2} \sum_{j,k=1}^{N_C}{\frac{d}{dt}f_\theta(t_i, x_j, v_k)}.\nonumber
\end{align}

\begin{figure}[ht]
    \centering
    \begin{subfigure}[t]{0.45\textwidth}
        \centering
        \includegraphics[width=\linewidth]{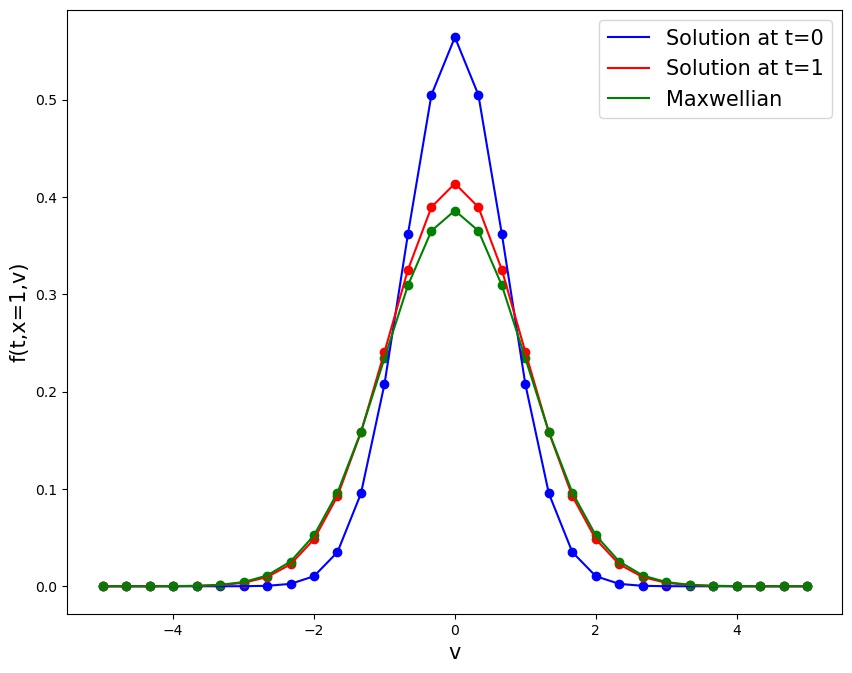} 
    \end{subfigure}
    %\hfill
    \begin{subfigure}[t]{0.45\textwidth}
        \centering
        \includegraphics[width=\linewidth]{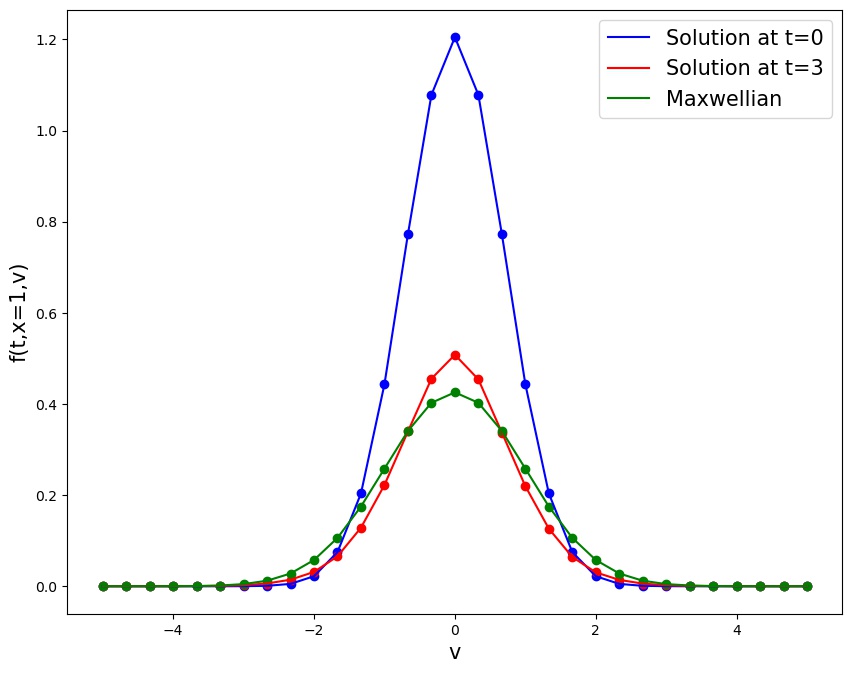} 
    \end{subfigure}
    \caption{Left: Relaxation to Maxwellian of f(t,x=1,v) from an initial condition $f_0^{(1)}$. Right: Relaxation to Maxwellian of f(t,x=1,v) from an initial condition $f_0^{(2)}$.}\label{fig:relaxation}
\end{figure}

For the kinetic Fokker--Planck equation, we consider two different initial conditions 
\begin{equation}
    \begin{split}
        \textbf{Test 1: }f_0^{(1)}(x,v) &= \frac{\exp(-v^2)}{\int_{-V}^V \exp(-v^2) dv}\nonumber\\
        \textbf{Test 2: }f_0^{(2)}(x,v) &= \frac{\cos(2\pi x)\exp(-v^2)}{\int_0^1\int_{-V}^V \cos(2\pi x)\exp(-v^2) dv},
    \end{split}
\end{equation}
where both initial conditions are normalized to obtain $\int f_0^{(i)} dxdv = 1$, for $i=1,2$. We show the relaxation to the global equilibrium for each initial condition in Figure \ref{fig:relaxation}.

\begin{figure}[ht]
\begin{center}
\centering
  \includegraphics[width=\textwidth]{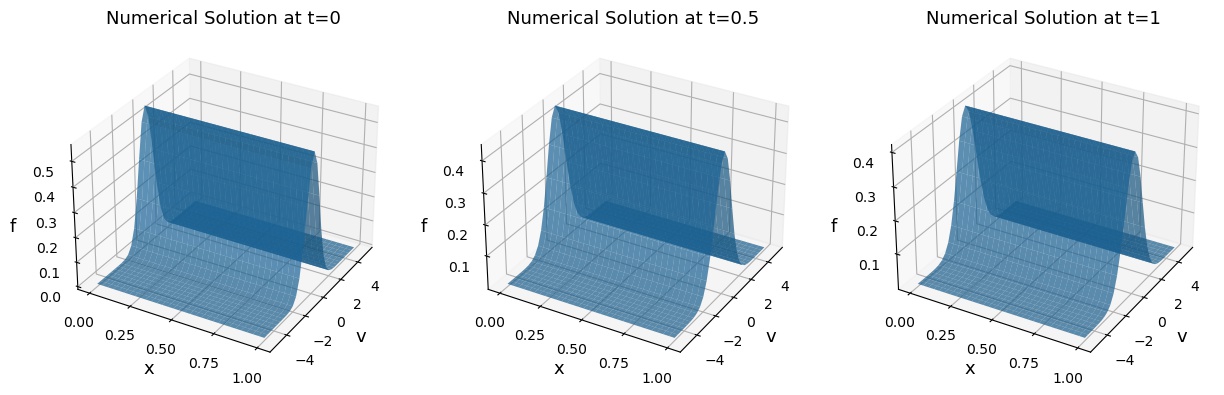}
  \caption{Numerical solutions of \textbf{Test 1} at $t=0, \frac{1}{2}, $ and $1.$}\label{fig:numerical_f1}
  \end{center}
\end{figure}

\textbf{Test 1.} In this example we take $f_0^{(1)}(x,v)$ as an initial condition with the periodic boundary condition. We set the time interval as $[0,T] = [0,1]$, thus the whole domain becomes $[0,1]\times[0,1]\times[-5,5]$. We set the diffusion coefficient by $q=1$, and the friction term by $p=1$. We show a numerical solution, computed by a method introduced in  \cite{wollman2008deterministic}, in Figure \ref{fig:numerical_f1}.

\begin{figure}[ht]
\begin{center}
\centering
  \includegraphics[height=0.4\textwidth,width=\textwidth]{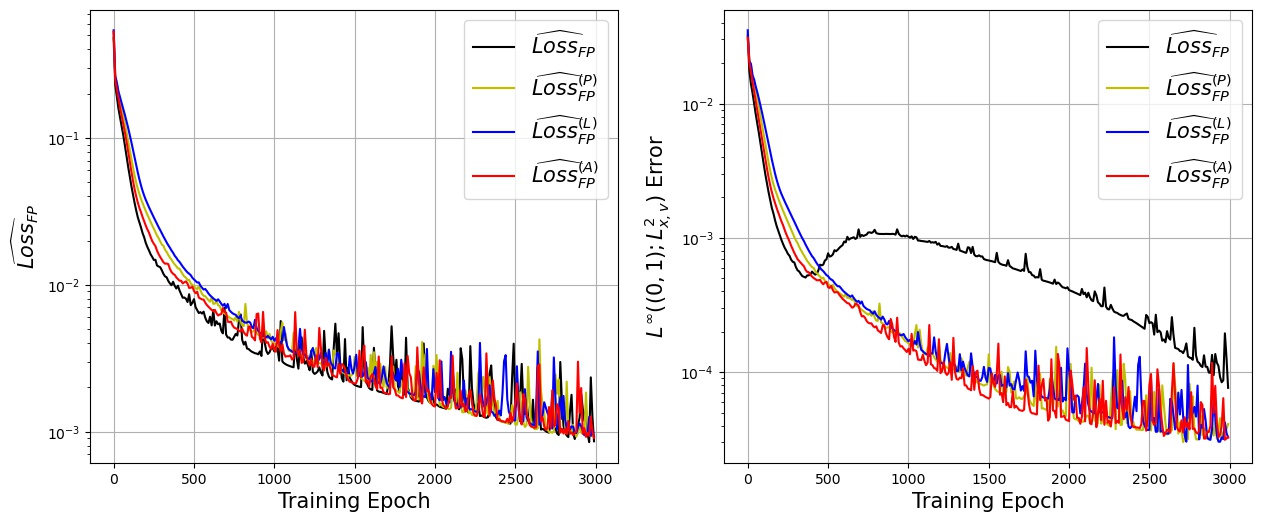}
  \caption{Left: Loss value in training epoch at log scale. Right: Actual error in training epoch at log scale.}\label{fig:FP_loss_error}
  \end{center}
\end{figure}

We train four identical neural networks with four different loss functions in \eqref{FP_loss}, $\widehat{Loss}_{FP}, \widehat{Loss}_{FP}^{(P)}, \widehat{Loss}_{FP}^{(L)}, \widehat{Loss}_{FP}^{(A)}$. We record the value of $\widehat{Loss}_{FP}$ for each training to see whether the constraints harm the original training problem or not. We also record the error of the neural network solution, and compare it to the numerical solution given in Figure \ref{fig:numerical_f1}. The results are summarized in the plots in Figure \ref{fig:FP_loss_error}. As we can see in the left panel, constraints have little to no effect when minimizing the loss value. On the other hand, the right panel shows that the proposed constrained optimization setting results in a far more accurate solution than the original unconstrained one as expected. One another interesting observation is that one can achieve the consistency of the loss and the error when using constrained loss function.

\begin{figure}[ht]
\begin{center}
\centering
  \includegraphics[height=0.4\textwidth,width=\textwidth]{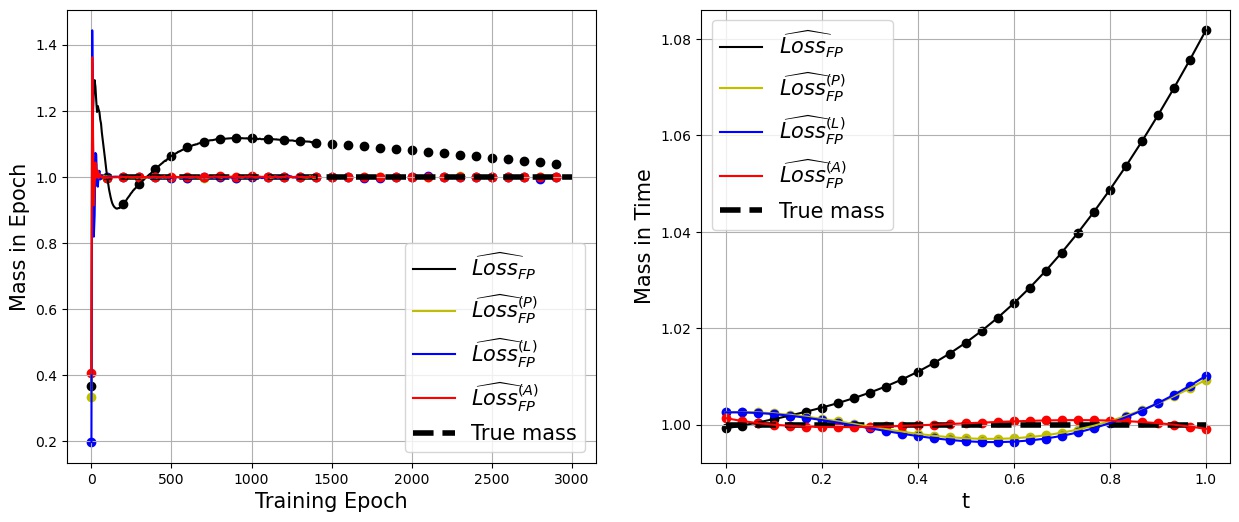}
  \caption{Left: Time averaged mass in training epoch. Right: Total mass in time after training is finished.}\label{fig:FP_mass}
  \end{center}
\end{figure}

Moreover, we show the mass conservation results in Figure \ref{fig:FP_mass}. Left panel shows the time-averaged mass \begin{equation}\nonumber
    \frac{1}{N_C} \sum_{i=1}^{N_C} \int_0^1\int_{-5}^5 f_\theta(t_i, x,v)dvdx,
\end{equation} for each training epoch. As we can see in the figure, three constrained settings give the accurate mass in a far earlier epoch than the unconstrained one. Right panel shows the mass in time $t$ \begin{equation}\nonumber
    \int_0^1\int_{-5}^5 f_\theta(t,x,v) dvdx, \text{ for $0\leq t\leq 1$,}
\end{equation} after the whole training process is done. Unconstrained one gives the increasing total mass in time while the constrained ones give almost constant mass in time. Overall performance considering both the accuracy and the mass conservation, is the best with $\widehat{Loss}_{FP}^{(A)}$, and the worst with $\widehat{Loss}_{FP}$.

\begin{figure}[ht]
\begin{center}
\centering
  \includegraphics[width=\textwidth]{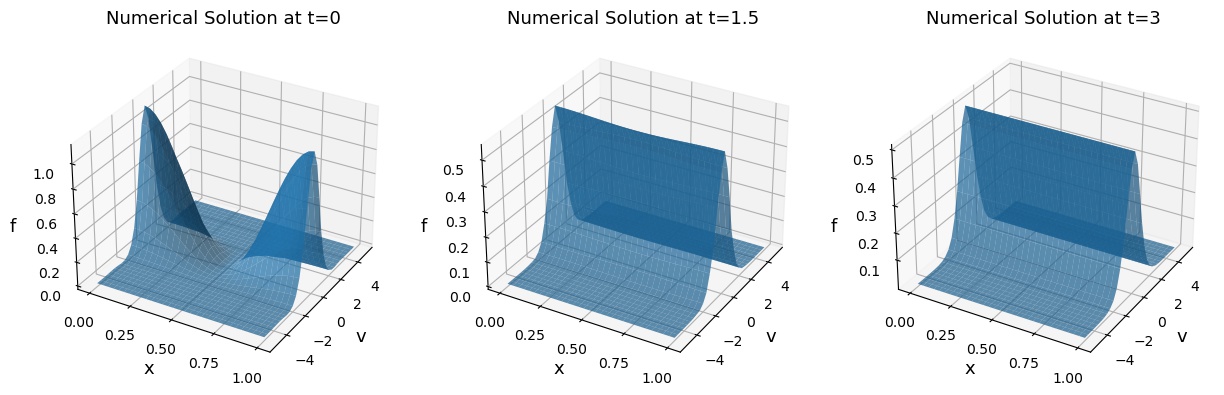}
  \caption{Numerical solutions of \textbf{Test 1} at $t=0, \frac{3}{2}, $ and $3.$}\label{fig:numerical_f2}
  \end{center}
\end{figure}

\textbf{Test 2.} Now we show the results of the test problem using $f_0^{(2)}(x,v)$ as an initial condition. We set the time interval as $[0,T] = [0,3]$ and the whole domain becomes $[0,3]\times[0,1]\times[-5,5]$. Both the diffusion and friction coefficients are set to $0.1$. We show a numerical solution, again computed by a method in \cite{wollman2008deterministic}, in Figure \ref{fig:numerical_f2}.

\begin{figure}[ht]
\begin{center}
\centering
  \includegraphics[height=0.4\textwidth,width=\textwidth]{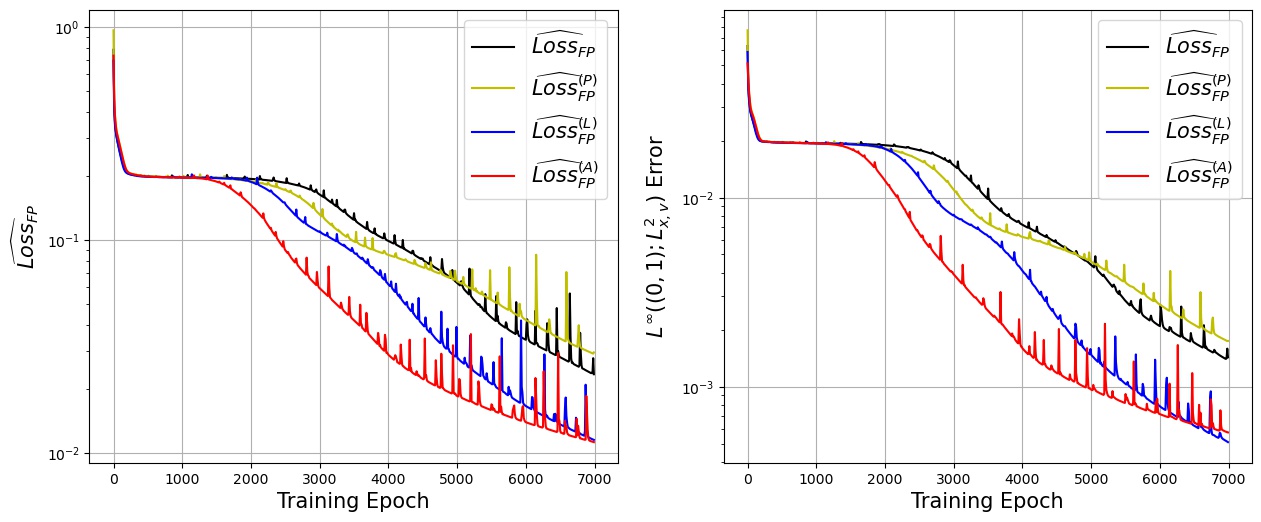}
  \caption{Left: Loss value in training epoch at log scale. Right: Actual error in training epoch at log scale.}\label{fig:FP_loss_error_f2}
  \end{center}
\end{figure}

\begin{figure}[ht]
\begin{center}
\centering
  \includegraphics[height=0.4\textwidth,width=\textwidth]{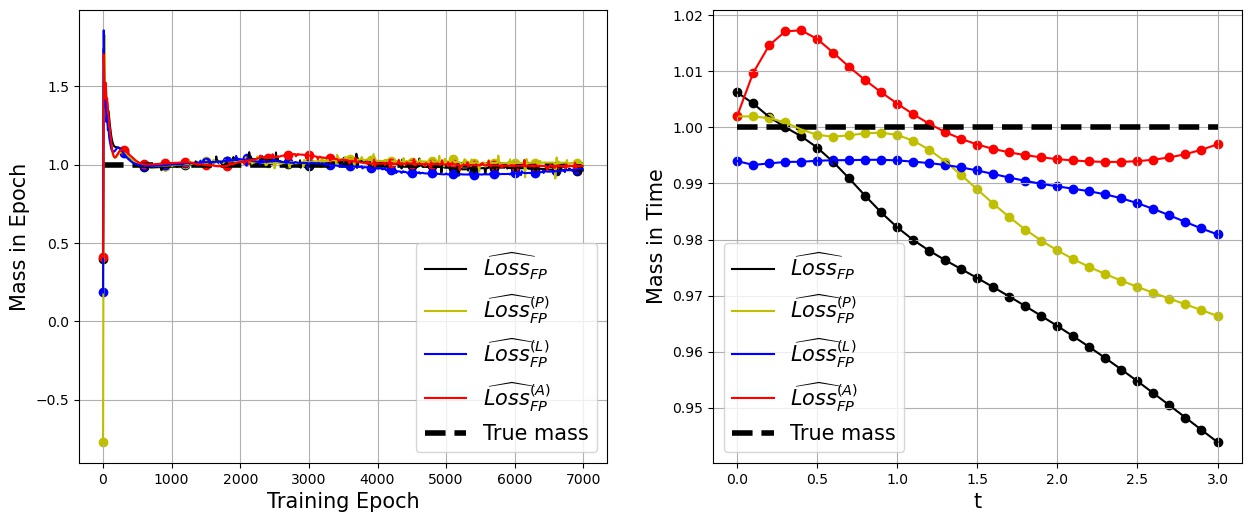}
  \caption{Left: Time averaged mass in training epoch. Right: Total mass in time after training is finished.}\label{fig:FP_mass_f2}
  \end{center}
\end{figure}

We present the same result for comparison of the loss $\widehat{Loss}_{FP}$ and the error $L^{\infty((0,3), L^2_{x,v})}$in Figure \ref{fig:FP_loss_error_f2}. In terms of error, we observe that the networks trained with $\widehat{Loss}_{FP}^{(L)}$, $\widehat{Loss}_{FP}^{(A)}$ are outperforming the others. Figure \ref{fig:FP_mass_f2} shows the result for the mass conservation, and still we observe a better performance of the proposed methods.

\subsection{Homogeneous Boltzmann equation}
In this subsection, we consider the homogeneous Boltzmann equation. We employ the test case considered in \cite{dimarco2018efficient, pareschi2000numerical, filbet2004accurate}. We consider the two dimensional Maxwellian molecules where $\sigma(|v-v_{*}|, w) = C$ and we set $C=1$.  We also truncate the velocity domain into $[-V,V]\times[-V,V]$ and we set the time interval as $[0,T]=[0,3]$. The equation reads 
\begin{equation}
    \begin{split}
        \partial_t f &= \frac{1}{\epsilon} Q(f,f), \\\nonumber
        f(0,v) &= \frac{v^2}{\pi} \exp(-v^2).
    \end{split}
\end{equation}
A well known analytic solution, called BKW solution, is given for all $t >0$ by \cite{bobylev1975exact, krook1977exact}:
\begin{equation} \label{BKW}
    f(t,v) = \frac{1}{2\pi S^2} \exp(\frac{-v^2}{2S}) (2S-1+\frac{1-S}{2S}v^2),
\end{equation}
where $S = 1-\frac{\exp(\frac{-t}{8})}{2}$.

As in the previous subsection, we sample the collocation points from each domain
\begin{equation}
    \begin{split}
        \{(t_i, v_{x,j}, v_{y,k})\}_{i,j,k=1}^{N_C} &\sim UNIF([0,3]\times[-5,5]\times[-5,5]) ,\nonumber\\
        \{(v_{x,j}, v_{y,k})\}_{j,k=1}^{N_I} &\sim UNIF([-5,5]\times[-5,5]),
    \end{split}
\end{equation}
where $N_C$, $N_I$ denote the number of sample points in the whole domain and initial domain, respectively. 
Then we discretize the loss in \eqref{loss_total_boltzmann}, and \eqref{proposed_boltzmann} as follows:
\begin{align}
        \widehat{Loss}_{B} &= \frac{4TV^2}{N_C^3} \sum_{i,j,k=1}^{N_C} (\partial_t f_\theta - Q(f,f))^2 \bigg|_{(t_i,v_{x,j}, v_{y,k})} \nonumber\\
        &+ \frac{4V^2}{N_I^2} \sum_{j,k=1}^{N_I} (f_\theta(0,v_{x,j},v_{y,k}) - f_0(v_{x_j}, v_{y,k}))^2, \nonumber\\
        \widehat{Loss}_{B}^{(P)} &= \widehat{Loss}_{B} + \sum_{l=1}^{4} \beta_l \frac{T}{N_C} \sum_{i=1}^{N_C}(c_l(t_i))^2, \label{Boltzmann_loss}\\
        \widehat{Loss}_{B}^{(L)} &= \widehat{Loss}_{B} + 
        \sum_{l=1}^{4} \sum_{i=1}^{N_C}\lambda(t_i)c_l(t_i), \nonumber\\
        \widehat{Loss}_{B}^{(A)} &= \widehat{Loss}_{B} + \sum_{l=1}^{4} \mu \frac{T}{N_C} \sum_{i=1}^{N_C}(c_l(t_i))^2  + \sum_{l=1}^{4} \sum_{i=1}^{N_C}\lambda(t_i)c_l(t_i), \nonumber
\end{align}
where
\begin{align}
       c_1(t;\theta) &= \frac{4V^2}{N_C^2}\sum_{j,k=1}^{N_C} \frac{d}{dt}f_\theta(t, v_{x_j}, v_{y,k}), \nonumber\\
       c_2(t;\theta) &= \frac{4V^2}{N_C^2}\sum_{j,k=1}^{N_C} \frac{d}{dt}f_\theta(t, v_{x_j}, v_{y,k})v_{x,j}, \nonumber\\
       c_3(t;\theta) &= \frac{4V^2}{N_C^2}\sum_{j,k=1}^{N_C} \frac{d}{dt}f_\theta(t, v_{x_j}, v_{y,k})v_{y,k}, \nonumber\\
       c_4(t;\theta) &= \frac{4V^2}{N_C^2}\sum_{j,k=1}^{N_C} \frac{d}{dt}f_\theta(t, v_{x_j}, v_{y,k})(v_{x,j}^2 + v_{y,k}^2),\nonumber
\end{align}
are the constraints.

\begin{figure}[ht]
\begin{center}
\centering
  \includegraphics[width=\textwidth]{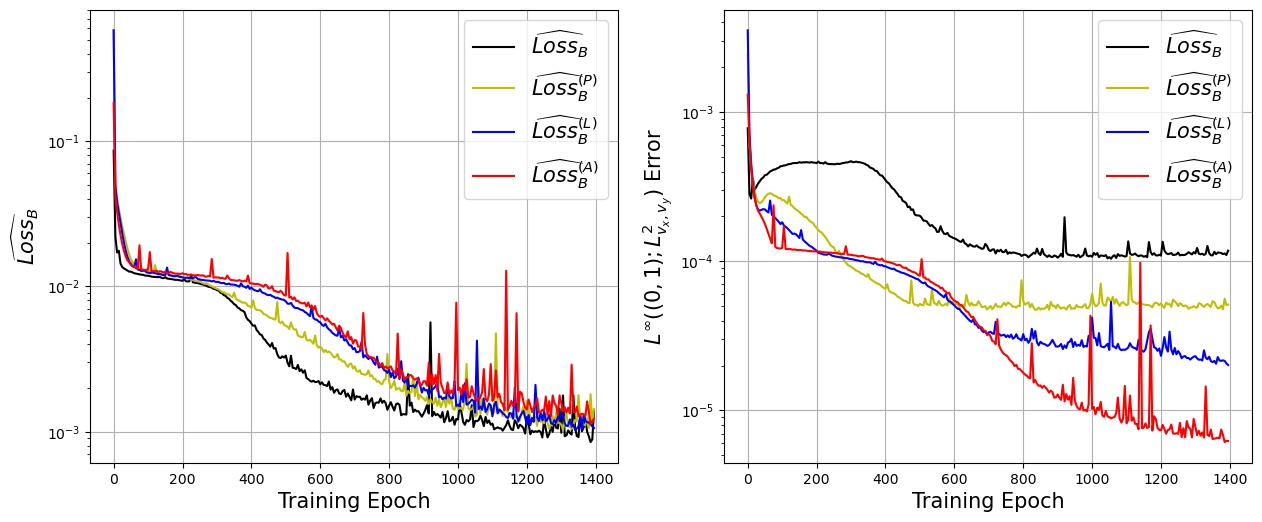}
  \caption{Left: Value of the loss $\widehat{Loss}_B$ in training epoch. Right: $L^{\infty}((0,1);L^2_{v_x,v_y})$ error in training epoch.}\label{fig:Boltzmann_loss}
  \end{center}
\end{figure}

For this problem, we train four neural networks with different loss functions in \eqref{Boltzmann_loss}. We summarize the value of loss in training epoch, and the error between the neural network solution and the BKW solution \eqref{BKW} in Figure \ref{fig:Boltzmann_loss}. Even though, the loss is smaller when we use the unconstrained one, the proposed conservative neural network significantly outperforms it in terms of the error $L^{\infty}((0,1);L^2_{v_x,v_y})$.

\begin{figure}[ht]
\begin{center}
\centering
  \includegraphics[width=\textwidth]{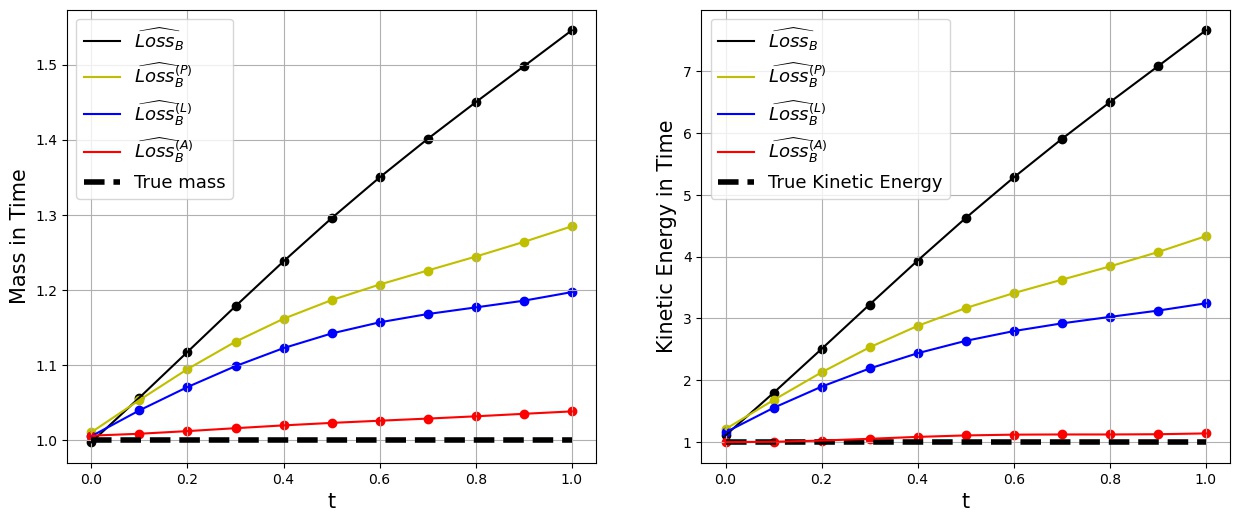}
  \caption{Left: Total mass in time after the training is finished. Right: Kinetic energy in time after the training is finished.}\label{fig:Boltzmann_mass_energy}
  \end{center}
\end{figure}

\begin{figure}[ht]
\begin{center}
\centering
  \includegraphics[width=\textwidth]{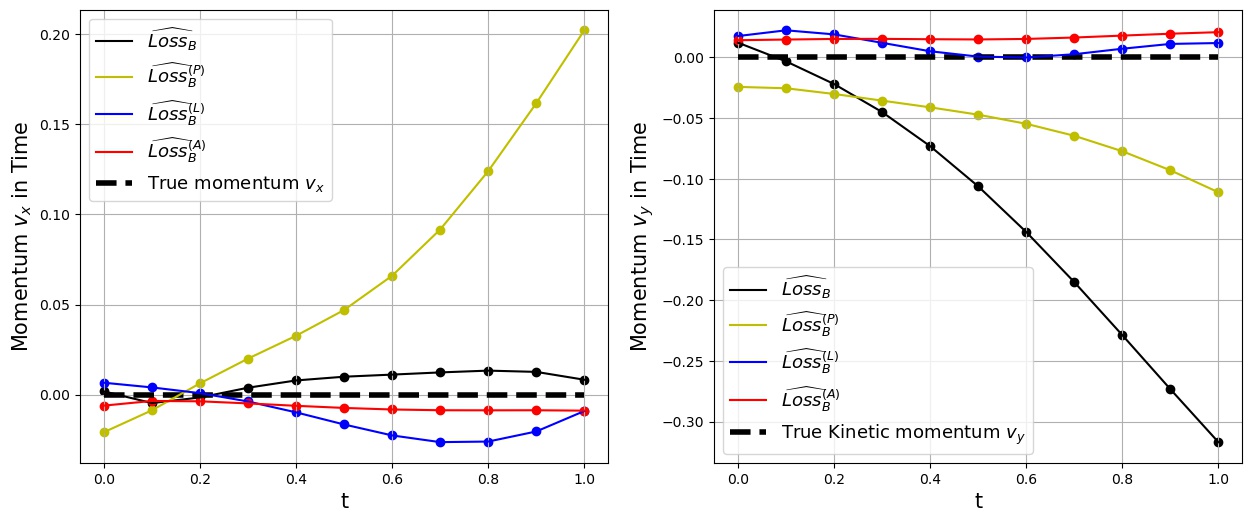}
  \caption{Momentum in time after the training is finished.}\label{fig:Boltzmann_momentum}
  \end{center}
\end{figure}

Regarding the conservation laws, we can see a superior performance of the proposed method. In Figure \ref{fig:Boltzmann_mass_energy}, we provide two plots which show the mass and the kinetic energy
\begin{equation}
    \text{Mass} : \int_{[-5,5]^2} f(t,v) dv, \quad \text{Kinetic Energy} : \frac{1}{2}\int_{[-5,5]^2} |v|^2 f(t,v) dv. \nonumber
\end{equation}
where both values are set to $1$ initially. Figure \ref{fig:Boltzmann_mass_energy} shows that the mass and the kinetic energy are only conserved in time when we train a neural network with $\widehat{Loss}_B^{(A)}$, which we proposed. Figure \ref{fig:Boltzmann_momentum} shows the conservation of momentum \begin{equation}
    \text{Momentum} : \int_{[-5,5]^2} vf(t,v) dv.\nonumber
\end{equation}
for neural networks trained with different loss functions in \eqref{Boltzmann_loss}. We also observe a better performance of the proposed loss functions in terms of the conservation of momentum. 

\section{Discussion}\label{discussion}
In this paper, we proposed a new framework for training a neural network for solving kinetic PDEs while satisfying conservative laws. We introduced a new class of loss functions based on the notion of constrained optimization where the constraints represent the physics conservation laws in the kinetic literature. The constraints are relaxed to an objective function, by the well known Lagrangian dual methods, so that we can optimize the objective function in a differentiable way. 

We validated our methodology through two equations, the kinetic Fokker--Plank equation and the homogeneous Boltzmann equation. For the kinetic Fokker--Plank equation, the only conserved quantity is mass, and the numerical results show that the proposed method gives far more accurate approximated solutions as well as the conservation property. For the 2-dimensional homogeneous Boltzmann equation, we impose four constrains, one for the mass, two for the momentum, and one for the kinetic energy. Again our methods significantly outperform the original unconstrained one in terms of both $L^{\infty}((0,T), L^2_{v})$-error and the conservation laws. 

The proposed methods are easy to implement and have almost the same computational cost as the original one, since the back-propagation of neural network dominates the cost. We also believe that this work can be extended to any kind of equations with conservation properties.

%For acknowledgements section, please don't number the section, please begin it with \section*{Acknowledgements}
\section*{Acknowledgments} This work was supported by the National Research Foundation of Korea (NRF) grant funded by the Korea government (MSIT) (NRF-2017R1E1A1A03070105, NRF-2019R1A5A1028324).

% You may incorporate your references as follows in your main tex file.
% Using BibTex is not recommended but can be handled.

\medskip
% The data information below will be filled by AIMS editorial staff
%Received xxxx 20xx; revised xxxx 20xx.
\medskip

\bibliographystyle{AIMS} %AMS references format with number references, and clickable URL links
\bibliography{bibliograph}

\end{document}